\documentclass[11pt]{article}
\usepackage{theorem}
\usepackage{amsmath}
\usepackage{graphicx}
\usepackage{amsfonts}
\usepackage{amssymb}
\newtheorem{proposition}{Proposition}
\newtheorem{theorem}[proposition]{Theorem}

{\theorembodyfont{\rmfamily}
\newtheorem{conjecture}[proposition]{Conjecture}
}
\newtheorem{lemma}[proposition]{Lemma}

\begin{document}

\author{Huai-Dong Cao and Bennett Chow}
\title{Recent Developments on the Ricci Flow}
\date{}
\maketitle
\begin{abstract}
This article reports recent developments of the research on Hamilton's Ricci flow and
its applications.
\end{abstract}

\section*{Introduction}

\footnotetext{\emph{1991 Mathematics Subject Classification\/.} Primary 58G11;
Secondary 53C21, 35K55.
\par
\noindent Authors partially supported by the NSF.} One of the fundamental 
problems in differential geometry is to find canonical
metrics on Riemannian manifolds. Here by canonical metrics we mean metrics
of constant curvature in various forms.  In turn, the existence of a canonical metric
on a manifold often has important topological implications.  A well-known example
is the uniformization theorem for closed surfaces. On the other hand, to find a
canonical metric on a given Riemannian manifold is often a very difficult problem.
For finding metrics of constant scalar curvature, it is the well-known Yamabe problem
(see R. Schoen's paper \cite{Schoen} and the references therein
for details). For metrics of constant 
Ricci curvature, i.e., Einstein metrics, one needs to solve Einstein's equation, which is 
extremely difficult in general.

In this article we shall describe the Ricci flow, or the parabolic Einstein equation,
introduced by Richard Hamilton in 1982 \cite{3-manifold} for producing Einstein metrics 
of positive scalar curvature and constant positive sectional curvature, and report on some 
of the more recent 
progress on the Ricci flow, especially the recent work of Hamilton which is focused on 
his program to understand Thurston's geometrization conjecture for three-manifolds using 
the Ricci flow.

In 1982, motivated by Eells and Sampson's work \cite{Eells-Sampson}
on the harmonic map flow in 1964, Hamilton introduced the Ricci flow
\begin{equation}
\frac\partial{\partial t}g_{ij}=-2R_{ij}.
\end{equation}
Notice that because of the minus sign in the front 
of the Ricci tensor in the equation, the solution metric to the Ricci flow shrinks in positive 
Ricci curvature direction while it expands in the negative Ricci curvature direction.
For example, on the $2$-sphere, any metric of positive Gaussian curvature will shrink to a 
point in finite time. Since the Ricci flow Eq.(1) does not preserve volume in
general, one often considers the \emph{normalized} Ricci flow defined
by
\begin{equation}
\frac\partial{\partial t}g_{ij}=-2R_{ij}+\frac2nrg_{ij},
\end{equation}
where $r=\left.  \int RdV\right/  \int dV$ is the average scalar curvature.
Under this normalized flow, which is equivalent to the (unnormalized) Ricci flow (1) by 
reparametrizing in time $t$
and scaling the metric in space by a function of $t$, the volume of the solution
metric is constant in time. Note also that Einstein metrics (i.e., $R_{ij}=c g_{ij}$) 
are fixed points of Eq.(2).

In \cite{3-manifold}, Hamilton showed that on a closed 
Riemannian $3$-manifold $M^{3}$ with initial metric of positive Ricci curvature, the
solution $g(t)$ to the normalized Ricci flow Eq.(2)
exists for all time and the metrics $g(t)$ converge exponentially fast, as time $t$ 
tends to the infinity, to a constant positive sectional curvature metric 
$g_{\infty}$ on $M^{3}.$ In particular, such a $M^3$ is necessarily diffeomorphic to a
quotient of the $3$-sphere by a finite group of isometries. It follows 
that given any homotopy $3$-sphere, if one can show that it admits 
a metric with positive Ricci curvature, then the Poincar\'{e} Conjecture would follow. 

Note the above result of Hamilton in particular implies that on a closed $3$-manifold with 
positive Ricci curvature, the curvature of the solution 
metric to the normalized Ricci flow is uniformly bounded in time $t$ for $0\leq t <\infty$.
However, on a general closed $3$-manifold, the solution of the normalized Ricci flow may
develop singularities,
meaning the curvature of the solution will become unbounded as time $t$ approach $T$, where 
$[0, T)$ is the maximal time interval for the existence of the normalized Ricci flow (2).
For example, if we take a dumbbell metric on the $3$-sphere $S^3$ with a neck like $S^2\times B^1$,
we expect the neck will shrink because the positive curvature in the $S^2$ direction
will dominate the slightly negative curvature in the $B^1$ direction. As a result, we expect the 
neck will pinch off in finite time, which means a singularity 
developing in the Ricci flow (topologically, neck pinching would correspond
to the prime decomposition of a $3$-manifold.)
 
Recall that Thurston's Geometrization Conjecture \cite{Thurston} states that every closed
$3$-manifold can be decomposed into pieces and each piece admits one of the 8 geometric 
structures, which would provide a
link between the geometry and topology of $3$-manifolds, analogous in spirit
to the case of surfaces. In particular, the conjecture has as a special case -
the Poincar\'{e} Conjecture.  

Hamilton's approach towards proving Thurston's Geometrization Conjecture is to study 
the Ricci flow on closed $3$-manifolds and try to show that given any initial Riemannian metric,
the Ricci flow evolves it to a geometric structure, after performing suitable
geometrical surgeries. More precisely, the program is to divide the study of the Ricci flow 
on closed $3$-manifolds into two parts. First, try to analyze 
singularities of the Ricci flow which develop in finite time well enough to enable one to
perform geometric surgeries before the singularities occur, which will
decompose the manifold, and then continue the solution (which will be nonsingular).  
Second, classify solutions to the normalized Ricci flow 
which exist for all time $t\in\lbrack0,\infty)$ and have uniformly bounded sectional curvature.
(These are called \emph{non-singular} solutions.)    

The rest of the paper is arranged as follows.  After reviewing some background
materials in section 1, the short time existence result
for the Ricci flow will be discussed in section 2. In section 3, we shall describe various 
convergence results of the normalized Ricci flow. 
In section 4, the recent
result of Hamilton on classification of non-singular solutions
to the normalized Ricci flow is stated. Section 5 describes the neck pinching
phenomenon for the Ricci flow on a closed $4$-manifold with non-negative isotropic curvature
where the situation is simpler. In section 6, we state the Harnack estimate for the Ricci flow on a positively 
curved manifold, an important tool in
analyzing singularities of the Ricci flow, and its important consequences.
Finally, in section 7, we discuss singularities of the Ricci flow for closed  $3$-manifolds
and explain how and why the study essentially reduces to the existence of a Harnack-type estimate
on a general closed $3$-manifold.

\section{Background material}

An important class of problems in Riemannian geometry is to understand the
interaction between the curvature and topology on a differentiable 
manifold.  A prime example of this interaction is the Gauss-Bonnet
formula on a closed surface $M^{2}$, which
says
\[
\int_{M}K\,dA=2\pi\,\chi\left(  M\right)  ,
\]
where $dA$ is the area element of a metric $g$ on $M$, $K$ is the Gaussian curvature
of $g$, and $\chi\left(  M\right)$ is the Euler characteristic of $M.$ 
 
To study the geometry of a differentiable manifold we
need an additional structure: the \emph{Riemannian metric}. The metric is
an inner product on each of the tangent spaces and tells
us how to measure angles and distances infinitesimally.  In local coordinates
$(x^1, x^2, \cdots, x^n)$, the metric $g$ is given by $\sum_{i,j}g_{ij}(x) dx^i\otimes
dx^j$, where $(g_{ij}(x))$ is a positive definite symmetric matrix at each point $x$.  
For a differentiable manifold one can differentiate functions. A Riemannian metric 
defines a natural way of
differentiating vector fields: \emph{covariant differentiation}. In Euclidean
space, one can change the order of differentiation. On a Riemannian manifold
the commutator of twice covariant differentiating vector fields is in general
nonzero and is called the \emph{Riemann curvature tensor}, which is a
$4$-tensor on the manifold.

For surfaces, the Riemann curvature tensor is equivalent
to the \emph{Gauss curvature} $K$, a scalar function.
In dimensions $3,$ or more, the Riemann curvature tensor $Rm$ is inherently a
tensor. In local coordinates, it is denoted by $R_{ijkl}$, which is anti-symmetric in 
$i$ and $k$ and in $j$ and $l$, and symmetric in the pairs $\{ij\}$ and $\{kl\}$.
Thus, it can be considered as a bilinear form on $2$-forms which is called the 
\emph{curvature operator}.
We now describe heuristically the various curvatures associated to the
Riemann curvature tensor. Given a point $x\in M^{n}$ and $2$-plane $\Pi$ in
the tangent space of $M$ at $x,$ we can define a surface $S$ in $M$ to be the
union of all geodesics passing through $x$ and tangent to $\Pi.$ In a
neighborhood of $x,$ $S$ is a smooth $2$-dimensional submanifold of $M.$ We
define \emph{the sectional curvature} $K\left(  \Pi\right)  $ of the $2$-plane
to be the Gauss curvature of $S$ at $x$:
\[
K\left(  \Pi\right)  =K_{S}\left(  x\right)  .
\]
Thus the sectional curvature $K$ of a Riemannian manifold associates to each
$2$-plane in a tangent space a real number. Given a line $L$ in a
tangent space, we can average the sectional curvatures of all planes through
$L$ to obtain the \emph{Ricci curvature} $Rc\left(  L\right)  .$ Likewise,
given a point $x\in M,$ we can average the Ricci curvatures of all lines in
the tangent space of $x$ to obtain the \emph{scalar curvature} $R(x).$
In local coordinates, the Ricci tensor is given by $R_{ik}=\sum_{jl} g^{jl}R_{ijkl}$
and the scalar curvature is given by $R=\sum_{ik} g^{ik}R_{ik}$, where $(g^{ij})$ 
is the inverse of the metric tensor $(g_{ij})$

Since the Ricci flow, which is the topic of this expository paper,
lies in the realm of parabolic partial differential equations, where the
prototype is the heat equation, below we give a brief review of the heat 
equation.

Let $(M^{n},g)$ be
a Riemannian manifold. Given a $C^{2}$ function $u:M\rightarrow\mathbb{R},$ 
its Laplacian
is defined in local coordinates $\left\{  x^{i}\right\}  $ to be
\[
\Delta u=\text{tr}_{g}\left(  \nabla^{2}u\right)  =g^{ij}\nabla_{i}\nabla
_{j}u,
\]
where  $\nabla
_{i}=\nabla_{\frac{\partial}{\partial x^{i}}}$ is its associated covariant
derivative (Levi-Civita connection.) 
We say that a $C^{2}$ function $u:M^{n}\times\lbrack0,T)\rightarrow
\mathbb{R},$ where $T\in(0,\infty],$ is a solution to the heat equation if
\[
\frac{\partial u}{\partial t}=\Delta u.
\]
One of the most
important properties satisfied by the heat equation is the maximum principle, which says that
for any smooth solution to the heat equation, whatever pointwise bounds hold
at $t=0$ also hold for $t>0.$

\begin{theorem}
\label{Max-prin1}\emph{(Maximum principle: pointwise bounds are preserved)}
Let $u:M^{n}\times\lbrack0,T)\rightarrow\mathbb{R}$ be a $C^{2}$ solution to
the heat equation on a complete Riemannian manifold. If $C_{1}\leq u\left(
x,0\right)  \leq C_{2}$ for all $x\in M,$ for some constants $C_{1},C_{2}%
\in\mathbb{R},$ then $C_{1}\leq u\left(  x,t\right)  \leq C_{2}$ for all $x\in
M$ and $t\in\lbrack0,T).$
\end{theorem}

This property, which exhibits the smoothing behavior of the heat equation,
follows from the following more general result.

\begin{proposition}
\label{max-prin2}Let $u:M^{n}\times\lbrack0,T)\rightarrow\mathbb{R}$ be a
$C^{2}$ function satisfying 
\[
\frac{\partial u}{\partial t}\leq\Delta u+X\cdot\nabla u,
\]
where the laplacian and dot product are defined with respect to a
time-dependent metric $g(t),$ and $X(t)$ is any time-dependent vector field.
If $u\left(  x,0\right)  \leq C$ for all $x\in M,$ for some constant
$C\in\mathbb{R},$ then $u\left(  x,t\right)  \leq C$ for all $x\in M$ and
$t\in\lbrack0,T).$
\end{proposition}

\emph{Idea of Proof. }At a point $(x,t)$ where $u$ attains its maximum at time
$t,$ we have $\Delta u\leq0$ and $\nabla u=0,$ which by the equation implies
$\partial u/\partial t\leq0.$

\bigskip

A lot of the results concerning the Ricci flow are proved using this version of
the maximum principle. An exception is when one obtains pointwise bounds for
the Ricci or sectional curvatures, where one uses a maximum principle for
systems, where the solutions (e.g., Ricci tensor or curvature operator) are no
longer functions but rather sections of a vector bundle. We refer the reader
to \cite{4-manifold} for details.

\section{The Short time existence for the Ricci flow}

Nonlinear heat equations first appeared in Riemannian geometry in 1964 when
Eells and Sampson \cite{Eells-Sampson} used the harmonic map heat flow to show that
any map between closed Riemannian manifolds, where the range has negative
sectional curvature, is homotopic to a harmonic map. It was the search for a
heat flow for Riemannian metrics which led Hamilton \cite{3-manifold} to
discover the Ricci flow in 1982. 
The basic idea of the Ricci flow is try to improve a given Riemannian metric by evolving
it by its Ricci curvature. 

Given a differentiable manifold $M,$ we
say that a one-parameter family of metrics $g\left(  t\right)  ,$ where
$t\in\lbrack0,T)$ for some $T>0,$ is a solution to the (unnormalized) Ricci flow if%
\[
\frac{\partial}{\partial t}g_{ij}=-2R_{ij}%
\]
at all $x\in M$ and $t\in\lbrack0,T).$ The minus sign in the equation makes 
the Ricci flow a \emph{forward} heat equation as we shall see below.
The factor $2$ is simply for normalization purpose.
 
Note that in local geodesic coordinates
$\{x^{i}\}$, we have
\[
g_{ij}(x)=\delta_{ij}-\frac{1}{3}R_{ipjq}x^{p}x^{q}+O\left(  \left|  x\right|
^{3}\right)  .
\]
Therefore%
\[
\Delta g_{ij}\left(  0\right)  =-\frac{1}{3}R_{ij}%
\]
where $\Delta$ is the standard euclidean laplacian.  Hence the 
Ricci flow is like 
the heat equation for a Riemannian metric
\[
\frac{\partial}{\partial t}g_{ij}=6\Delta g_{ij}.
\]

The practical study of the Ricci flow is made possible by the following
short-time existence result.

\begin{proposition}
Given any smooth compact Riemannian manifold $(M,g_{o}),$ there exists a
unique smooth solution $g(t)$ to the Ricci flow defined on some time interval
$t\in\lbrack0,\epsilon)$ such that $g(0)=g_{o}.$
\end{proposition}

We remark that the Ricci flow is a weakly parabolic system where degeneracy comes from 
the gauge invariance
of the equation under diffeomorphisms. Therefore, short time existence does 
not follow from general theory.
Richard Hamilton's original proof of the short time existence was involved and used the Nash-Moser inverse
function theorem. Soon after, D. DeTurck \cite{DeTurck} substantially
simplified the short-time existence proof by breaking the diffeomorphism
invariance of the equation (which causes difficulty in directly applying
standard theory to prove short-time existence.)

On the other hand, when $M$ is a complex manifold and the initial metric
$g_{0}$ is K\"ahler the Ricci flow is strictly parabolic. This is due to the
fact that the gauge group of biholomorphisms is a much smaller group compared
with the full diffeomorphism group. In fact, in the K\"ahler case the Ricci
flow can be reduced to a strictly parabolic scalar equation of Monge-Amp\`ere
type (see \cite{Ba} and \cite{Ca1}). Hence the short time existence follows easily.

\section{Convergence results}

Given that short-time existence holds for any smooth initial metric, one of
the main problems concerning the Ricci flow is to determine under what
conditions the solution to the normalized equation exists for all time and 
converges to a constant curvature
metric. Results in this direction have been
established under various curvature assumptions, most of them being some sort
of positive curvature. Since the Ricci flow (1) does not preserve volume in
general, one often considers, as we mentioned in the introduction, the  
normalized Ricci flow (2):
\[
\frac\partial{\partial t}g_{ij}=-2R_{ij}+\frac2nrg_{ij}.
\]
Under this flow,  the volume of the solution $g(t)$ is independent of time.

To study the long-time existence of the normalized Ricci flow, it is important
to know what kind of curvature conditions are preserved under the equation.
In general, the Ricci flow tends to preserve some kind of positivity of
curvatures. For example, positive scalar curvature is preserved in all
dimensions. This follows from applying the maximum principle, Proposition
\ref{max-prin2}, to the evolution equation for scalar curvature $R,$ which is
\[
\frac{\partial}{\partial t}R=\Delta R+2\left|  R_{ij}\right|^{2}.
\]
In dimension $3$, positive Ricci curvature is preserved under the Ricci flow.
This is a special feature of dimension $3$ and is related to the fact that the
Riemann curvature tensor may be recovered algebraically from the Ricci tensor
and the metric in dimension $3$. Positivity of sectional curvature is not
preserved in general. However, the stronger condition of positive
curvature operator is preserved under the Ricci flow. Recall that the
Riemann curvature tensor may be considered as a self-adjoint map
$Rm:\wedge^{2}M\rightarrow\wedge^{2}M.$ We say that a metric $g$ has positive
(non-negative) curvature operator if the eigenvalues of $Rm$ are positive
(non-negative.) We remark that positivity of curvature operator implies the
positivity of the sectional curvature (and in dimension $3,$ the two conditions
are equivalent.) Finally in the K\"{a}hler case, the condition of positive
holomorphic bisectional curvature, which is a weaker condition than positive
sectional curvature, is also preserved (see \cite{Ba} and \cite{Mo}).

Although the condition of positive scalar curvature is preserved in all dimensions, no
convergence results are known for metrics satisfying this condition 
except in dimension $2$. In dimension $3$ one expects necks to pinch off corresponding
to the prime decomposition of a $3$-manifold. 

To illustrate what the Ricci flow can do, we first state Hamilton's 1988 result on 
surfaces \cite{surface}. It is interesting to note that this result was 
obtained a few years after his celebrated work on $3$-manifolds (see below) 
and the proof is, in some sense,  even harder than the latter one.

\begin{theorem}
Let $M$ be a closed surface. Then for any initial metric $g_{0}$ on $M$, the
solution to the normalized Ricci flow exists for all time. Moreover,

\begin{enumerate}
\item  If the Euler characteristic of $M$ is non-positive, then the solution
metric $g(t)$ converges to a constant curvature metric as $t\rightarrow\infty.$

\item  If the scalar curvature $R$ of the initial metric $g_{0}$ is positive,
then the solution metric $g(t)$ converges to a positive constant curvature
metric as $t\rightarrow\infty.$
\end{enumerate}
\end{theorem}

For surfaces with non-positive Euler characteristic,
the proof was based primarily on maximum principle estimates for the scalar
curvature. The proof for the case of metrics with positive
scalar curvature  is highly non-trivial and uses what 
are known as Entropy and Harnack estimates (both
depend on having $R>0$). Convergence in the remaining case of a compact 
surface with positive Euler characteristic and initial metric with scalar curvature
changing sign was proved in \cite{2-sphere} by extending the techniques in
\cite{surface}. Bartz, Struwe, and Ye \cite{BSY} and Hamilton
\cite{isoperimetric} have both given new proofs of this result; the former
based on the Aleksandrov reflection method and the latter based on an
isoperimetric estimate. 

We now turn our attention to Hamilton's famous work on $3$-manifold in 1982 
\cite{3-manifold}.

\begin{theorem}
\emph{(Positive Ricci topological spherical space form)} Let $(M^{3},g_{0})$
be a closed Riemannian $3$-manifold with positive Ricci curvature. Then there
exists a unique solution to the normalized Ricci flow $g(t)$ with $g(0)=g_{0}$
for all time and the metrics $g(t)$ converge exponentially fast to a constant
positive sectional curvature metric $g_{\infty}$ on $M^{3}.$ In particular,
$M^{3}$ is diffeomorphic to a spherical space form.
\end{theorem}

As a consequence, such a $3$-manifold $M$ is necessarily diffeomorphic to a
quotient of the $3$-sphere by a finite group of isometries. It follows 
that given any homotopy $3$-sphere, if one can show that it admits 
a metric with positive Ricci curvature, then the Poincar\'{e} Conjecture would follow. 
More generally, the Elliptization Conjecture would follow from showing that any
closed $3$-manifold with finite fundamental group admits a metric with
positive Ricci curvature.

Prior to Hamilton's work, it was not known that a  constant
positive sectional curvature metric (or equivalently, Einstein metric
in dimension 3) exists on a closed $3$-manifold with positive Ricci curvature.
(Of course, this would follow from the Poincar\'{e} conjecture
and Spherical Space Form Conjecture.)
In fact, it is very important and difficult to show the existence of an Einstein 
metric on a given manifold in general. Yau's solution of Calabi conjecture gives a 
powerful method for finding K\"ahler-Einstein metrics of non-positive scalar curvature, 
while Hamilton's Ricci 
flow produces Einstein metrics of positive scalar curvature.

Remarks on the proof of the theorem: The proof involves obtaining certain
strong a priori estimates for the curvature and its derivatives. To give the
reader a taste of the type of estimates involved, we outline the main estimate
for the Ricci curvatures below. Recall that a metric is Einstein if the Ricci
tensor is proportional to the metric:
\[
R_{ij}=\frac{1}{n}R g_{ij}.
\]
When $n\geq3,$ this implies the scalar curvature $R$ is constant. It is a
special feature of dimension $3$ that Einstein implies constant sectional
curvature. In dimension 3, one of the main estimates is the following, which says that the
metric is close to Einstein at points where the scalar curvature is large.

\begin{lemma}
There exist constants $\varepsilon>0$ and $C<\infty$ depending only on the
initial metric such that%
\[
\frac{\left|  R_{ij}-\frac{1}{3}R g_{ij}\right|  ^{2}}{R^{2}}\leq C\,R^{-\varepsilon}%
\]
\end{lemma}

We note that the left-hand side is scale-invariant (i.e., does not change when the metric
is multiplied by a positive constant.) The right-hand side is small when $R$ is large.
Thus, if we can prove $R_{\min}\rightarrow\infty,$ then we would have a
scale-invariant pointwise measure of the difference of the metric from one having
constant sectional curvature tends to zero uniformly: $\left.  \left|
R_{ij}-\frac{1}{3}R g_{ij}\right|  ^{2}\right/  R^{2}\rightarrow0.$ To show that
$R_{\min}\rightarrow\infty$ requires more estimates, and we refer the reader
to \cite{3-manifold} for details. In addition, higher derivative estimates for
the curvatures are required to prove convergence in $C^{\infty}.$ 

Recall that in all dimensions, the condition of positive curvature operator is
preserved under the Ricci flow. In 1986 Hamilton \cite{4-manifold} proved the following

\begin{theorem}
If $(M,g_{o})$ is a compact $4$-manifold with positive curvature operator,
then there exists a unique solution $g(t)$ with positive curvature operator to
the normalized Ricci flow with $g(0)=g_{o}$ for all time $t\in\lbrack
0,\infty)$ such that as $t\rightarrow\infty$ the metric $g(t)$ converges in
$C^{\infty}$ to a smooth metric $g_{\infty}$ on $M$ with constant positive
sectional curvature.
\end{theorem}

Again, the significance of this result is that the Ricci flow produced a constant positive
sectional curvature metric on a compact $4$-manifold with positive curvature operator.
As a consequence, such a  $4$-manifold must be diffeomorphic to either $S^{4}$ 
or $\mathbb{R}P^{4}.$
Hamilton has conjectured that the above long-time existence and convergence
result holds in all dimensions. This is still unresolved for $n\geq5,$
although when the initial metric has sufficiently pointwise pinched positive
sectional curvatures, convergence results have been obtained by G. Huisken
\cite{Huisken}, C. Margerin \cite{Margerin}, and S. Nishikawa \cite{Nishikawa}. 
These results generalized the well-known differentiable sphere theorems; see \cite{AM}
for a history. In
dimension $4,$ Hamilton's result has been generalized by H.\ Chen
\cite{2-positive} to the case of $2$-positive curvature operator, which means the sum
of every two eigenvalues of the curvature operator is positive.

To a large extent, one should view the surface case as a special case of the Ricci flow on
K\"ahler manifolds, since Riemann surfaces can be regarded as 1-dimensional
complex manifolds. For convergence of the Ricci flow on compact K\"ahler
manifolds, the first author \cite{Ca1} proved the following result in 1985.

\begin{theorem}
Let $M$ be a compact K\"{a}hler manifold with definite first Chern class
$c_{1}(M)$. If $c_{1}(M)=0$, then for any initial K\"{a}hler metric $g_{0}$,
the solution to the normalized Ricci flow exists for all time and converges to
a Ricci flat metric as $t\rightarrow\infty.$ If $c_{1}(M)<0$ and the initial
metric $g_{0}$ is chosen to represent the negative of the first Chern class,
then the solution to the normalized Ricci flow exists for all time and
converges to an Einstein metric of negative scalar curvature as $t\rightarrow
\infty.$
If $c_{1}(M)>0$ and the initial metric $g_{0}$ is chosen to represent the
first Chern class, then the solution to the normalized Ricci flow exists for
all time.
\end{theorem}

The proof is partly based on the a priori estimates of Yau for Monge-Amp\`{e}%
re equations. In case of $c_{1}(M)>0$, the convergence is open even when $M$
has positive holomorphic bisectional curvature. For more recent progress on this 
problem, see \cite{Ca3} and \cite{CH}.

For the study of Ricci flow on noncompact K\"{a}hler manifolds, W.-X Shi has done very nice works.
In \cite{Sh2}, Shi used the Ricci flow to prove that a complete non-compact
K\"{a}hler manifold $X^{n}$ with bounded and positive holomorphic bisectional
curvature must be biholomorphic to $\mathbb{C}^{n},$ provided the manifold has
Euclidean volume growth and the average scalar curvature has quadratic decay.
This provides a partial affirmative answer to a conjecture of Yau. In
\cite{Sh3}, he dropped the condition of Euclidean volume growth and assumed
that the average scalar curvature decays like $1/R^{1+\epsilon}$ and concluded
that the manifold $X^{n}$ is biholomorphic to a pseudoconvex domain in
$\mathbb{C}^{n}.$

\section{Non-singular solutions on $3$-manifolds}

As we have seen in the last section, the Ricci flow had important topological consequence
for closed $3$-manifolds of positive Ricci curvature. Now we turn our attention to 
the study of the Ricci flow on general closed $3$-manifolds.  

Although long-time existence and convergence hold for the Ricci flow on
compact $2$-dimensional manifolds with arbitrary initial metrics, this is
certainly not the case in dimensions $3$ and higher, where singularities
develop in finite time for certain initial metrics. Here, by singularities of
the Ricci flow we mean solutions to the normalized Ricci flow (2) with unbounded
curvature. A typical example is when the manifold, say $S^n$,
is shaped like a dumbbell.
That is, there are two larger spherical regions to the left and right joined
by a cylindrical region (the `neck') in the middle. The topology of the neck
is that of a cylinder $S^{n-1}\times [-1,1]$. When $n\geq 3$, the sphere
$S^{n-1}$ has positive  intrinsic curvature, which makes the neck shrink 
(provided the neck doesn't open up too fast,)
leading to a singularity.
Therefore the study of the Ricci flow on general $3$-manifolds is much more 
complicated and difficult.

However, one can divide the study of the Ricci flow on $3$-manifolds
into two parts.
First, try to analyze the singularities which develop in finite time well 
enough to enable one to
perform geometric surgeries before the singularities occur, which will
decompose the manifold, and then continue the solution. Second, classify solutions 
to the normalized Ricci flow which exist for all
time $t\in\lbrack0,\infty)$ and have uniformly bounded sectional curvature
(these are called \emph{non-singular} solutions.). In dimension $3$ much
progress has been made by Hamilton in both directions, especially the second part 
of the program, where Hamilton \cite{non-singular} recently proved that non-singular
solutions have geometric decompositions, which we shall describe below.

A solution $(M,g)$ to the Ricci flow (normalized or unnormalized) on a time
interval $[0,T)$ is said to be \emph{maximal} if it cannot be extended past
time $T.$ In this case either

\begin{enumerate}
\item $T=\infty$ and $\sup_{M\times\lbrack0,T)}\left|  Rm\right|  <\infty$, or

\item $\sup_{M\times\lbrack0,T)}\left|  Rm\right|  =\infty.$
\end{enumerate}

In the first case we say that the solution is \emph{non-singular, }and in the
second case we say that the solution is \emph{singular }and that a singularity
occurs at time $T\leq\infty.$

\subsection{Examples}

By the results stated in the previous sections, the normalized Ricci flow on a
closed manifold is non-singular when the initial metric is

\begin{enumerate}
\item [N1.] any metric on a surface

\item [N2.]  a metric on a $3$-manifold with positive Ricci curvature

\item [N3.] a locally homogeneous metric on a $3$-manifold (see Isenberg-Jackson
\cite{Isenberg-Jackson})

\item [N4.] a metric on a $4$-manifold with positive curvature operator

\item [N5.] a metric on an $n$-manifold with sufficiently pointwise-pinched
positive sectional curvatures

\item  [N6.] a K\"{a}hler metric on a compact complex $n$-manifold with first Chern class
 $c_{1}=0$ or $c_{1}<0.$
\end{enumerate}

In each case, we actually have convergence to an Einstein
metric as $t\rightarrow T.$ In contrast, the unnormalized Ricci flow on a
compact manifold is singular when the initial metric is

\begin{enumerate}
\item [U1.] any metric on a surface with $\chi>0$

\item[U2.] a metric on a $n$-manifold with positive scalar curvature

\item[U3.] a locally homogeneous metric on a $3$-manifold of class $SU(2)$ or
$S^{2}\times\mathbb{R}$

\item[U4.] a K\"{a}hler metric on a compact complex $n$-manifold with $c_{1}>0.$
\end{enumerate}

Note that in cases U1 and U3 the singularity is \emph{removable } by normalizing
the flow. The same is true in case U2 if
 the initial metric is as in either
case
N2,
N4 or N5.

\subsection{Non-singular solutions}

In dimension 3, non-singular solutions are now well-understood 
topologically
through Hamilton's work \cite{non-singular}.
In particular, closed $3$-manifolds which admit a non-singular solution can
also be decomposed into geometric pieces. This focuses the study of the Ricci
flow on closed $3$-manifolds on the analysis of singularities. Now we present
the main result in Hamilton's work \cite{non-singular}.

\begin{theorem}
If a closed differentiable $3$-manifold $M^{3}$ admits a non-singular solution
$g(t)$ to the normalized Ricci flow, i.e., a solution which exists for all
time and has uniformly bounded sectional curvature, then $M^{3}$ has a
decomposition into geometric pieces. In particular, $M^{3}$ is diffeomorphic
to either

\begin{enumerate}
\item [1.] a Seifert fibered space

\item [2.] a spherical space form $S^{3}/\Gamma$

\item [3a.] a flat manifold

\item [3b.] a hyperbolic (constant negative sectional curvature) manifold

\item [4.] the union along incompressible tori of finite volume hyperbolic
manifolds and Seifert fibered spaces.
\end{enumerate}
\end{theorem}

Hamilton actually proves the following. If $(M^{3},g(t))$ is a non-singular
solution to the normalized Ricci flow on a closed $3$-manifold, then exactly
one of the following occurs:

\begin{enumerate}
\item (Sequential collapse) There exists a sequence of times $t_{i}%
\rightarrow\infty$ such that the metrics $g(t_{i})$ collapse, i.e.,
\[
\lim_{i\rightarrow\infty}\rho\left(  t_{i}\right)  ^{2}\max_{M\times\{t_{i}%
\}}\left|  Rm\right|  =0,
\]
where $\rho\left(  t_{i}\right)  =\max_{x\in M}$inj$\left(  x,t_{i}\right)  $
is the maximum injectivity radius of any point for the metric $g(t_{i}).$ By
Cheeger-Gromov theory \cite{CG}, $M^{3}$ admits an $F$-structure 
and is topologically a
graph manifold. This implies $M^{3}$ is a Seifert fibered space. (See \cite{A1}
for an exposition of Cheeger-Gromov theory.)

\item (Exponential convergence to a spherical space form) The solution $g(t)$
converges exponentially fast in every $C^{k}$-norm to a constant positive
sectional curvature metric $g_{\infty}$ on $M^{3}.$ This implies $M^{3}$ is
diffeomorphic to a space form $S^{3}/\Gamma.$

\item (Sequential convergence to a flat or hyperbolic manifold) There exists a
sequence of times $t_{i}\rightarrow\infty$ and self-diffeomorphisms $\phi_{i}$
of $M^{3}$ such that the sequence of pulled back metrics $\phi_{i}^{\ast
}g(t_{i})$ converge to a constant (either zero or negative) sectional
curvature metric $g_{\infty}$ on $M^{3}$ as $i\rightarrow\infty.$

\item (Toroidal decomposition into hyperbolic and Seifert fibered pieces)
There exist a finite collection of complete noncompact $3$-manifolds $\left\{
H_{\alpha}^{3},h_{\alpha}\right\}  $ with constant negative sectional
curvature and finite volume, and smooth $1$-parameter families of
diffeomorphisms (into) $\psi_{\alpha}(t):H_{\alpha}^{3}\rightarrow M^{3}$ such
that the pulled-back metrics $\psi_{\alpha}(t)^{\ast}g(t)$ converge to
$h_{\alpha}$ as $t\rightarrow\infty.$ Moreover, $M^{3}$ can be decomposed into
two time-dependent $3$-manifolds $M_{1}(t)$ and $M_{2}(t)$ ($M_{j}(t)$ and
$M_{j}(t^{\prime})$ are isotopic for all $t$ and $t^{\prime},$ $j=1,2$) whose
intersection is their mutual boundary, which consists of a finite disjoint
union of $2$-tori. Geometrically, the metrics pulled-back from $M_{1}(t)$
converge to the hyperbolic pieces, i.e., $\left(  \psi_{\alpha}^{-1}M_{1}%
,\psi_{\alpha}^{\ast}g\right)  (t)$ converge to $\cup_{\alpha}\left(
H_{\alpha}^{3},h_{\alpha}\right)  $ as $t\rightarrow\infty,$ while $\left(
M_{2}(t),g(t)\right)  $ collapses. Topologically, the boundary tori are all
incompressible, and all of the homomorphisms induced by the inclusion maps
$i_{\ast}:\pi_{1}\left(  N\right)  \rightarrow\pi_{1}\left(  P\right)  ,$
where $N=M_{1}(t)\cap M_{2}(t),$ $M_{1}(t)$ or $M_{2}(t)$ and $P=M_{1}(t),$
$M_{2}(t)$ or $M,$ are injective.
\end{enumerate}

The most difficult case is the last one where hyperbolic limits occur. To show
that $\pi_{1}(H_{i})$ injects into $\pi_{1}(M)$ for each hyperbolic limit,
Hamilton uses a minimal surface argument, part of which is reminiscent of
Schoen and Yau's \cite{Schoen-Yau} use of minimal surfaces to classify compact
$3$-manifolds with positive scalar curvature up to homotopy type (and modulo
the spherical space-form conjecture.) The idea is to represent a non-trivial
element of the kernel by a minimal disk and show that the rate of change of
the area of the minimal disk under the Ricci flow is uniformly bounded from
above by a negative constant leading to a contradiction since the initial disk 
has finite area. The main estimate is
to show that the length of the boundary of the minimal disk bounds the area of
the disk. Hamilton also proves that $M$ may be written as the union of two
manifolds (depending on time) $M_{1}(t)$ and $M_{2}(t)$ whose intersection is
their mutual boundary which is a finite disjoint union of tori such that
$\left(  \psi_{i}(t)^{-1}\left[  M_{1}(t)\right]  ,\psi_{i}(t)^{\ast
}g(t)\right)  $ converges to $\cup_{i=1}^{k}\left(  H_{i},h_{i}\right)  $ as
$t\rightarrow\infty,$ and $\left(  M_{2}(t),g(t)\right)  $ collapses as
$t\rightarrow\infty.$

\section{Necks in a $4$-dimensional case}

The first part of Hamilton's program, the analysis of singularities, is more
difficult than the second. However, in dimension $4$ some aspects of it are
simpler than in dimension $3.$ Recently Hamilton \cite{PIC} has classified
compact $4$-manifolds with non-negative\emph{\ isotropic curvature} using the
Ricci flow where singularities (pinching necks) develop in finite time. A
manifold of dimension $n\geq4$ has non-negative isotropic curvature if and
only if for every orthonormal $4$-frame $\left\{  e_{i}\right\}  _{i=1}^{4},$
we have
\[
R_{1313}+R_{1414}+R_{2323}+R_{2424}\geq2R_{1234},
\]
where $R_{ijkl}=\left\langle R(e_{i},e_{j})e_{l},e_{k}\right\rangle .$ In 1988
Micallef and Moore \cite{MM} showed that a simply-connected compact $n$-manifold with
positive isotropic curvature is homeomorphic to $S^{n}$ by studying the index
of minimal $2$-spheres. In dimension $4,$ Hamilton relaxed the condition that
$M$ be simply-connected and also obtained a classification up to
diffeomorphism. In particular, he proved

\begin{theorem}
Let $M$ be a compact $4$-manifold with no essential incompressible space-form,
that is, with the property that if $N$ is a submanifold of $M$ diffeomorphic
to a space-form $S^{3}/\Gamma$ such that $i_{\ast}:\pi_{1}(N)\rightarrow
\pi_{1}(M)$ is an injection, then $\Gamma$ is isomorphic to either $\left\{
1\right\}  $ or $\mathbb{Z}_{2}.$ If $M$ admits a metric with non-negative
isotropic curvature, then $M$ is diffeomorphic to either $S^{4}$ or the
connected sum of copies of $\mathbb{R}P^{4},$ $S^{3}\times S^{1},$ and
$S^{3}\widetilde{\times}S^{1}$ (the unique non-orientable $S^{3}$ bundle over
$S^{1}.$) The converse is also true.
\end{theorem}

The proof of this result uses a difficult analysis of the system of PDEs
satisfied by the Riemann curvature operator to show that the only
singularities which can develop are necks pinching off. Here a neck is a
(usually small) piece of the manifold which metrically is close to a quotient 
by isometries of a constant
multiple of the standard metric on a long thin finite cylinder $S^{3}\times
\lbrack-L,L],$ where $L\gg1$ (one also needs to consider $\mathbb{Z}_{2}$
quotients.) At a suitable time before a singularity forms, Hamilton
geometrically performs a surgery at the neck which either is of the form
$S^{3}\times B^{1}\rightarrow B^{4}\times S^{0}$ or removes a $\mathbb{R}%
P^{4}$ summand and then discards any component of the new manifold which is
diffeomorphic to either $S^{4},$ $\mathbb{R}P^{4},$ $S^{3}\times S^{1},$ or
$S^{3}\widetilde{\times}S^{1}.$ He then continues the Ricci flow starting with
the new $4$-manifold and metric and shows that after only a finite number of
surgeries one obtains the empty set. One of the basic ingredients in analyzing
the singularities and performing the surgeries is the Harnack estimate.

\section{The Harnack estimate}

Harnack inequalities are fundamental in the study of elliptic and parabolic
partial differential equations. In the Ricci flow it is basic in the
understanding of singularities. In \cite{Harnack} Hamilton proved a complicated
\emph{differential} Harnack inequality for the Ricci flow on $n$-manifolds
using the method of Li and Yau developed in \cite{Li-Yau} in 1986.

\begin{theorem}
\emph{(The Harnack estimate for the Ricci flow)}
If $\left( M,g(t)\right) $ is a solution to the Ricci flow on a compact
manifold with non-negative operator, then for any vector field $W$ and
$2$-form $U,$ we have 
\[
Z=M_{ij}W^iW^j+2P_{ijk}U^{ij}W^k+R_{ijkl}U^{ij}U^{kl}\geq 0, 
\]
where 
\[
P_{ijk}=\nabla _iR_{jk}-\nabla _jR_{ik} 
\]
and 
\[
M_{ij}=\Delta R_{ij}-\frac 12\nabla _i\nabla
_jR+2g^{kp}g^{lq}R_{ikjl}R_{pq}-g^{kl}R_{ik}R_{jl}+\frac 1{2t}R_{ij}. 
\]
\end{theorem}

In the K\"{a}hler case, a similar Harnack estimate was proved by the first
author \cite{Ca2} under the weaker curvature assumption of non-negative
holomorphic bisectional curvature.  See \cite{geometric} for a geometric 
interpretation of Hamilton's Harnack estimate.
As we shall see below, the Harnack estimate is important because it is the major tool
in analyzing singularities of the Ricci flow. 

The first important consequence of the Harnack estimate is to relate singularities 
of the Ricci flow to Ricci solitons,  which are special solutions of the Ricci 
flow moving along the equation by diffeomorphisms, that is, 
where the metric $g(t)$ is
the pull-back of the initial metric $g(0)$ by a $1$-parameter family of 
diffeomorphism $\phi\left( t\right)$ generated by a vector field on manifold $M$.
(In the PDE
literature, they are more commonly referred to as self-similar solutions.)
It turns out that on a Ricci soliton,  
the Harnack-Li-Yau quadratic $Z$ in the above theorem vanishes.  
Using the Harnack estimate, one can show that under certain conditions, limits of dilations 
(i.e., blow up) of
singularities  of the Ricci flow are (gradient) Ricci solitons (see \cite{soliton} and \cite{Ca3}).
Hence, the study of Ricci solitons also becomes important.

The equations that characterize a Ricci soliton metric are given by 
$\nabla_{i}W_{j}=\nabla_{j}W_{i}=R_{ij},$ where $W=\{W^i\}$ is a vector field on $M$.
If $W$ is a gradient vector field, then we have a \emph{gradient} Ricci soliton.
The first example of a \emph{steady} gradient Ricci soliton was found
by Hamilton \cite{surface} in dimension $2$. It is defined on $\mathbb{R}^{2}$
and has the form
\[
ds^{2}=\frac{dx^{2}+dy^{2}}{1+x^{2}+y^{2}}.
\]
This Ricci soliton is called the \emph{cigar soliton} $\Sigma,$ because it is
asymptotic to a flat cylinder at infinity and has maximal curvature at the
origin. Higher dimensional examples of rotational symmetric K\"{a}hler-Ricci
solitons have been found recently by the first author (see \cite{Ca3} and
\cite{Ca4}).

Another important consequence of the Harnack estimate is the following 
\emph{Little Loop Lemma} for $3$-manifolds with non-negative sectional 
curvature (see \cite{survey}, section 15). Roughly speaking, it states that 
there is no little geodesic loop in a large flat region. More precisely, we have

\begin{lemma}(\emph{Little Loop Lemma})
Let $(M,g(t))$ be a solution to the Ricci flow on a compact $3$-manifold with
non-negative sectional curvature. There exists a universal constant $A>0$ and
a constant $B>0$ depending only on the initial metric $g_{0}$ such that if $x$
is a point where $R(y,t)\leq\frac{A}{\rho^{2}}$ for all $y\in B_{\rho}(x)$ at
some time $t$ and for some $\rho,$ then we have the injectivity radius estimate
\[
\text{inj}\left(  x\right)  \geq B\rho,
\]
at that time $t.$
\end{lemma}

As we shall see in the next section, the Little Loop Lemma can be used to rule 
out certain types of singularities.

Hamilton conjectures that this result holds in all dimensions under the
assumption of non-negative curvature operator.\footnote{The proof in
\cite{survey} for all dimensions is apparently not complete. However, Hamilton
has announced a complete proof when $n=3.$}

\section{Singularities in dimension $3$}

One reason that the Harnack estimate is relevant, even under the restrictive
assumption of non-negative curvature operator, is that in dimension $3,$
limits arising from dilating singularities have non-negative sectional
curvature (which is the same as non-negative curvature operator when $n=3.$)
In particular the Harnack estimate is one of the estimates used in proving
the following classification of $3$-dimensional singularities given by
Hamilton in \cite{survey}.

\begin{proposition}
If $\left(  M,g(t)\right)  $ is a solution to the Ricci flow on a compact
$3$-manifold where a singularity develops in finite time $T,$ then there
exists a sequence of dilations of the solution which converges to a quotient
by isometries of either $S^{3},$ $S^{2}\times\mathbb{R},$ or $\Sigma
\times\mathbb{R},$ where $\Sigma$ is the cigar soliton.
\end{proposition}

Hamilton also conjectured that the limits which are quotients of $\Sigma
\times\mathbb{R}$ cannot occur and outlined an approach for proving this. If
this is true, then there always exists a sequence of dilations which converge to a
quotient either $S^{3}$ or $S^{2}\times\mathbb{R}.$ In the first case, the
manifold has to be a topological space form. In the second case, one sees a
neck, on which one wants to perform a surgery. One would then like to show
that after a finite number of surgeries, the solution becomes a non-singular solution.

To see how one might be able to rule out limits which are quotients of $\Sigma
\times\mathbb{R}$, notice that at its infinity, $\Sigma$ is asymptotic to a cylinder, 
hence the limits which
are quotients of $\Sigma\times\mathbb{R}$ have short (geodesic) loops in relatively flat
regions and so does the manifold right before taking the limits. 
Now for a solution to the Ricci flow on a compact $3$-manifold which has
non-negative sectional curvature, these short geodesic loops can be ruled out by 
the Little Loop Lemma stated in the last section. The problem is that even though 
the limits of dilations have non-negative sectional curvature, the solution metric right 
before taking the limit does not have non-negative sectional curvature, hence the Harnack estimate
and little loop lemma do not apply directly. 
Thus we need to have the Little Loop Lemma for solutions 
to the Ricci flow
with possible small negative curvature somewhere, as in the case right before taking 
the limit of dilations.

In conclusion, the limits
$\Sigma\times\mathbb{R}$ and $\Sigma\times S^{1}$ can be ruled out in general
if the following conjecture of Hamilton is true.

\begin{conjecture}
A Harnack-type estimate holds for the Ricci flow on compact $3$-manifolds with
arbitrary initial metric.
\end{conjecture}

Thus proving this conjecture is a crucial step in Hamilton's program to 
understand Thurston's geometrization conjecture. One reason for believing 
that this conjecture may be true is the
following estimate saying that in a sense the
sectional curvatures tend to positive (see \cite{survey} or \cite{Ivey}.)

\begin{proposition}
Let $(M,g(t))$ be a solution to the Ricci flow on a compact $3$-manifold.
There exists constants $C<\infty$ and $c>0$ depending only on the initial
metric $g_{0}$ such that if at some point $(x,t)$, the scalar curvature
 $R(x,t)\geq C,$  then
the minimum sectional curvature $K_{\min}$  satisfies
\[
K_{\min}\geq-c\frac{R}{\ln R}%
\]
at $(x,t).$
\end{proposition}

Thus, if the scalar curvature $R$ is large at a point, then $\left|  K_{\min}\right|  $ is much
smaller than $R$ at that point.

\bigskip

There are also many other important works on the Ricci flow which we have not
discussed in this article, including the works by Bemelmans-Min-Oo-Ruh, R.
Bryant, Cafora-Isenberg-Jackson, T. Ivey, D. Knopf, N. Koiso,
Leviton-Rubinstein, Min-Oo, Y. Shen, L.-F. Wu, D. Yang, and R. Ye, to name a
few. Finally we mention that another approach to the Geometrization 
Conjecture has been studied by M. Anderson \cite{A2}.

\textbf{Acknowledgment}. We are especially grateful to Richard Hamilton for
continually teaching us about the Ricci flow over the many years, and to our
advisor Shing-Tung Yau for getting us started in the subject and keeping us going.

\noindent Department of Mathematics, Texas A\&M University, College Station,
TX 77843. \emph{E-mail address}: cao@math.tamu.edu\newline 

\noindent Department of Mathematics, University of Minnesota, Minneapolis, MN
55455. \emph{E-mail address}: bchow@math.umn.edu
\end{document}